\documentclass[12pt,a4paper,oneside,reqno]{article}

\usepackage[utf8]{inputenc}
\usepackage[english]{babel}
\usepackage[left=1in,right=1in,top=1in,bottom=1in]{geometry}
\usepackage{amsfonts,amssymb,amsmath,amsthm}
\usepackage{xcolor}
\usepackage[normalem]{ulem}
\usepackage{cite}
\usepackage{hyperref}
\hypersetup{colorlinks=true,linkcolor=blue,citecolor=blue,urlcolor=blue}
\usepackage{fancyhdr}
\newtheorem{theorem}{Theorem}
\newtheorem{lemma}[theorem]{Lemma}

\newcommand{\doi}[1]{\href{https://doi.org/#1}{doi:#1}}
\newcommand{\arxiv}[1]{\href{https://arxiv.org/abs/#1}{arXiv:#1}}

\begin{document}

\noindent\fbox{\parbox{\dimexpr\textwidth-2\fboxsep-2\fboxrule\relax}{\small
\textbf{Corrected version verified in Lean.}  This is a corrected version of the paper
I.~Kaygorodov, \emph{$\delta$-Superderivations of semisimple finite-dimensional Jordan
superalgebras}, Math. Notes 91 (2012), no.~2, 187--197, \doi{10.1134/S0001434612010208}.}}

\bigskip

\begin{center}
{\Large $\delta$-Superderivations of semisimple\\[2pt] finite-dimensional Jordan superalgebras}

\bigskip

{\bf Ivan Kaygorodov}

\end{center}

\medskip

\noindent{\small\textbf{Abstract.} In the paper, a complete description of the $\delta$-derivations and
the $\delta$-superderivations of semisimple finite-dimensional Jordan superalgebras over an algebraically
closed field of characteristic $p\neq2$ is given. In particular, new examples of
nontrivial $\frac{1}{2}$-derivations and odd $\frac{1}{2}$-superderivations are given that are not operators of
right multiplication by an element of the superalgebra.}

\medskip

\noindent{\small\textbf{Keywords:}  
Jordan superalgebra, $\delta$-derivation,
$\delta$-superderivation, $\frac{1}{2}$-derivation, algebraically closed field.}

\bigskip

\noindent {\bf MSC2020:} 17C70 (primary); 17A36, 17A70, 17C20, 17C55 (secondary).

\bigskip

\setcounter{tocdepth}{1}
\tableofcontents

\section*{Introduction}
\addcontentsline{toc}{section}{Introduction}

Antiderivations, i.e., linear mappings $\mu$ of an algebra $A$ such that
\[
\mu(xy)=-\mu(x)y-x\mu(y),
\]
were studied in~\cite{1,2}. The notion of $\delta$-derivation, i.e., a linear mapping
$\varphi$ of an algebra $A$ such that
\[
\varphi(xy)=\delta(\varphi(x)y+x\varphi(y))
\]
for a chosen element $\delta$ of the ground field and for arbitrary elements $x,y\in A$, was later
introduced in papers of Filippov. He considered $\delta$-derivations of prime Lie, alternative, and
Maltsev non-Lie algebras~\cite{3,4,5}. After this, $\delta$-derivations of Jordan algebras and
superalgebras and also of simple Lie superalgebras were considered by Kaygorodov~\cite{6,7,8,9},
$\delta$-superderivations of superalgebras of the Jordan bracket were treated by Kaygorodov and
Zhelyabin~\cite{10}, and Zusmanovich considered $\delta$-superderivations of prime Lie
superalgebras~\cite{11}. Note that a detailed survey devoted to the study of $\delta$-derivations and
$\delta$-superderivations is given in~\cite{12}.

In the present paper, a complete description of the $\delta$-derivations and $\delta$-superderivations of
simple finite-dimensional nonunital Jordan superalgebras over an algebraically closed field of
characteristic $p>2$ is given. As a result, we have a complete description of the $\delta$-derivations
and $\delta$-superderivations of the semisimple finite-dimensional Jordan superalgebras over an
algebraically closed field of characteristic distinct from~$2$. In particular, it is proved that the semisimple finite-dimensional Jordan superalgebras
over an algebraically closed field of characteristic zero have no nontrivial $\delta$-derivations and
$\delta$-superderivations.

\section{Main facts and definitions}

Let $F$ be a field of  characteristic $p\neq2$. An algebra $A$ over a field $F$ is said to be
Jordan if it satisfies the identities
\[
xy=yx,\qquad (x^2y)x=x^2(yx).
\]
Let $G$ be the Grassmann algebra over $F$ given by the generators $1,\xi_1,\dots,\xi_n,\dots$ and the
defining relations $\xi_i^2=0$ and $\xi_i\xi_j=-\xi_j\xi_i$. The elements
\[
1,\quad \xi_{i_1}\xi_{i_2}\cdots\xi_{i_k},\quad i_1<i_2<\dots<i_k,
\]
form a basis of the algebra $G$ over $F$. Denote by $G_0$ and $G_1$ the subspaces spanned by the products
of even and odd lengths, respectively; then $G$ can be represented as the direct sum of these subspaces,
$G=G_0\oplus G_1$. Here the relations $G_iG_j\subseteq G_{i+j\,(\mathrm{mod}\,2)}$, $i,j=0,1$, hold. In
other words, $G$ is a $\mathbb Z_2$-graded algebra (or a superalgebra) over~$F$.

Suppose now that $A=A_0\oplus A_1$ is an arbitrary superalgebra over $F$. Consider the tensor product
$G\otimes A$ of $F$-algebras. The subalgebra
\[
G(A)=G_0\otimes A_0+G_1\otimes A_1
\]
of $G\otimes A$ is referred to as the Grassmann envelope of the superalgebra~$A$.

Let $\Omega$ be a variety of algebras over $F$. A superalgebra $A=A_0\oplus A_1$ is referred to as an
$\Omega$-superalgebra if its Grassmann envelope $G(A)$ is an algebra in $\Omega$. In particular,
$A=A_0\oplus A_1$ is referred to as a Jordan superalgebra if its Grassmann envelope $G(A)$ is a Jordan
algebra. Further, for a homogeneous element $x$ of a superalgebra $A=A_0\oplus A_1$, we set $p(x)=i$ if
$x\in A_i$ and denote by $x_i$ the projection of an element $x\in A=A_0\oplus A_1$ to~$A_i$.

A classification of the finite-dimensional Jordan superalgebras over an algebraically closed field of
characteristic zero was given in the papers by Kantor and Kac~\cite{13,14}. Later on, Racine and
Zel'manov~\cite{15} described the simple finite-dimensional Jordan superalgebras with semisimple even
part over an algebraically closed field of characteristic distinct from~$2$. In the paper by Martinez
and Zelmanov~\cite{16}, a classification of the simple unital finite-dimensional Jordan superalgebras
with nonsemisimple even part over an algebraically closed field of characteristic $p>2$ was given. The
simple nonunital finite-dimensional Jordan superalgebras over an algebraically closed field of
characteristic $p>2$ were described by Zelmanov~\cite{17}. The structure of semisimple finite-dimensional
Jordan superalgebras over an algebraically closed field of characteristic $p\neq2$ was also presented
in that paper.

Let us present some examples of Jordan superalgebras.

\subsection*{1.1. Superalgebra of vector type $J(\Gamma,D)$}

Let $\Gamma=\Gamma_0\oplus\Gamma_1$ be an associative supercommutative unital superalgebra with nonzero
even derivation $D$. Consider the direct sum $J(\Gamma,D)=\Gamma\oplus\Gamma x$ of spaces, where
$\Gamma x$ stands for an isomorphic counterpart of the space $\Gamma$. Then the multiplication operation
on $J(\Gamma,D)$, which we denote by $\cdot$, is defined by the formulas
\begin{center}
$a\cdot b=ab,\quad 
a\cdot bx=(ab)x,\quad 
ax\cdot b=(-1)^{p(b)}(ab)x,\quad
ax\cdot bx=(-1)^{p(b)}(D(a)b-aD(b)),$
\end{center}
where $a$ and $b$ are homogeneous elements in $\Gamma$ and $ab$ stands for the product in~$\Gamma$.

In what follows, we are interested in the case of $\Gamma=B(m)$, where
\[
B(m)=F[a_1,\dots,a_m\mid a_i^p=0]
\]
stands for the algebra of truncated polynomials in $m$ even variables over the field of characteristic
$p>2$.

\subsection*{1.2. Kac superalgebra}

The simple nine-dimensional Kac superalgebra $K_9$ over a field of characteristic~$3$ is defined as
follows:
\begin{center}
$K_9=A\oplus M,\quad (K_9)_0=A,\quad (K_9)_1=M,$ $A=Fe+Fuz+Fuw+Fvz+Fvw,\quad M=Fz+Fw+Fu+Fv.$
\end{center}
The multiplication is defined by the following conditions:
\[
e^2=e,\ \ e\text{ is the unit element of the algebra }A,\quad em=\tfrac12m\ \text{ for any } m\in M,
\]
\[
[u,z]=uz,\quad [u,w]=uw,\quad [v,z]=vz,\quad [v,w]=vw,
\]
\[
[z,w]=e,\quad [u,z]w=-u,\quad [v,z]w=-v,\quad [u,z][v,w]=2e,
\]
and the other nonzero products are obtained from those presented above either by applying one of the
skew symmetries $z\leftrightarrow w$ and $u\leftrightarrow v$ or by applying the simultaneous
substitution $z\leftrightarrow u$, $w\leftrightarrow v$. Note that the superalgebra $K_9$ is not unital.

\subsection*{1.3. Superalgebra $V_{\frac{1}{2}}(Z,D)$}

Let $Z$ be an associative and commutative $F$-algebra with unit $e$ and with a derivation
$D\colon Z\to Z$ satisfying the following two conditions:
\begin{itemize}
\item[i)] $Z$ has no proper nonzero $D$-invariant ideals;
\item[ii)] $D$ annihilates elements of the form $Fe$ only.
\end{itemize}
Consider $Zx$ as an isomorphic counterpart of the algebra $Z$. Equip the vector space
\[
V(Z,D)=Z+Zx
\]
with the structure of a superalgebra. Let $A=Z$ and $M=Zx$ be the even and odd parts, respectively.
Define the product $\cdot$ as follows:
\[
a\cdot b=ab,\qquad a\cdot bx=\tfrac12(ab)x,\qquad ax\cdot bx=D(a)b-aD(b)
\]
for arbitrary elements $a,b\in Z$. Denote by $V_{\frac{1}{2}}(Z,D)$ the superalgebra thus obtained.

As was noted above, for a chosen element $\delta$ of the ground field, by a $\delta$-derivation of a
superalgebra $A$ we mean a linear mapping $\varphi\colon A\to A$ satisfying the condition
\[
\varphi(xy)=\delta(\varphi(x)y+x\varphi(y))
\]
for arbitrary $x,y\in A$.

By the centroid $\Gamma(A)$ of the superalgebra $A$ we mean the set of all linear mappings
$\chi\colon A\to A$ satisfying the condition
\[
\chi(ab)=\chi(a)b=a\chi(b)
\]
for arbitrary $a,b\in A$.

The definition of $1$-derivation coincides with the standard definition of derivation; a
$0$-derivation is an arbitrary endomorphism $\varphi$ of the algebra $A$ such that $\varphi(A^2)=0$. It
is clear that every element of the centroid of a superalgebra is a $\frac{1}{2}$-derivation.

A nonzero $\delta$-derivation $\varphi$ is said to be a nontrivial $\delta$-derivation if
$\delta\neq0,1$ and $\varphi\notin\Gamma(A)$.

By a superspace we mean a $\mathbb Z_2$-graded space. Considering the space of endomorphisms
$\mathrm{End}(A)$ of the superalgebra $A$, we can define a $\mathbb Z_2$-grading by assuming that the even
mappings are the mappings $\varphi$ such that $\varphi(A_i)\subseteq A_i$ and the odd ones are the
mappings $\varphi$ such that $\varphi(A_i)\subseteq A_{i+1}$. A homogeneous element $\psi$ of the
superspace $\mathrm{End}(A)$ is referred to as a superderivation if
\[
\psi(xy)=\psi(x)y+(-1)^{p(x)p(\psi)}x\psi(y)
\]
for any homogeneous elements $x,y\in A_0\cup A_1$.

For a chosen element $\delta\in F$, introduce the notion of $\delta$-superderivation of a superalgebra
$A=A_0\oplus A_1$. A homogeneous linear mapping $\varphi\colon A\to A$ is referred to as a
$\delta$-superderivation if
\[
\varphi(xy)=\delta(\varphi(x)y+(-1)^{p(x)p(\varphi)}x\varphi(y))
\]
for any homogeneous elements $x,y\in A_0\cup A_1$.

Consider a Lie superalgebra $A=A_0\oplus A_1$ with multiplication $[\,\cdot\,,\cdot\,]$ and choose an
element $x\in A_i$. Then $R_x\colon y\mapsto[x,y]$ is a superderivation of the superalgebra $A$, and its
parity is $p(R_x)=i$.

By the supercentroid $\Gamma_s(A)$ of the superalgebra $A$ we mean the set of all homogeneous linear
mappings $\chi\colon A\to A$ satisfying the condition
\[
\chi(ab)=\chi(a)b=(-1)^{p(a)p(\chi)}a\chi(b)
\]
for arbitrary homogeneous elements $a$ and~$b$.

The definition of $1$-superderivation coincides with the definition of ordinary superderivation; a
$0$-superderivation is an arbitrary endomorphism $\varphi$ of the superalgebra $A$ such that
$\varphi(A^2)=0$. A nonzero $\delta$-superderivation is said to be nontrivial if $\delta\neq0,1$ and
$\varphi\notin\Gamma_s(A)$.

It is easy to see that an even $\delta$-superderivation is a $\delta$-derivation. Below we use this fact
without mentioning it explicitly.

\section{$\delta$-derivations and $\delta$-superderivations of simple finite-dimensional Jordan
superalgebras}

In this section, we consider $\delta$-derivations and $\delta$-superderivations of nonunital simple
finite-dimensional Jordan superalgebras over an algebraically closed field of characteristic $p>2$,
i.e., of the superalgebras $K_9$ over a field of characteristic~$3$ and of the superalgebra
$V_{\frac{1}{2}}(Z,D)$ over a field of characteristic $p>2$. Finally, with regard to results of~\cite{6}
and~\cite{8}, we shall have a complete description of the $\delta$-derivations and the
$\delta$-superderivations of simple finite-dimensional Jordan superalgebras over an algebraically closed
field of characteristic distinct from~$2$.

Below, denote by $\langle x,y\rangle$ the vector space spanned by vectors $x$ and~$y$.

\begin{lemma}\label{l1}
The superalgebra $K_9$ has no nontrivial $\delta$-derivations and $\delta$-superderivations.
\end{lemma}

\begin{proof}
Let $\varphi$ be a nontrivial $\delta$-derivation of the superalgebra $K_9$. It is easy to see that
\[
0=\varphi(me-em)=2\delta m\varphi(e)_1
\]
for $m\in(K_9)_1$. Hence, since $\varphi(e)_1u=0$ and $\varphi(e)_1v=0$, it follows that
$\varphi(e)_1=0$.

Consider the restriction of $\varphi$ to $(K_9)_0$. Note that $(K_9)_0$ is a Jordan algebra of a
bilinear form; as was proved in~\cite{6},
\[
\varphi(x)_0=\alpha x,\quad x\in(K_9)_0,\ \alpha\in F,
\]
for $\delta=\frac{1}{2}$. In particular, $\varphi(e)_0=\alpha e$. If $\delta\neq\frac{1}{2}$, then
$\varphi((K_9)_0)_0=0$.

If $\delta\neq\frac{1}{2}$ and $m\in(K_9)_1$, then
\[
\varphi(m)=2\varphi(em)=2\delta e\varphi(m)=2\delta\varphi(m)_0+\delta\varphi(m)_1.
\]
This gives $\varphi(m)=0$ and $\varphi(m_1m_2)=\delta\varphi(m_1)m_2+\delta m_1\varphi(m_2)=0$ for
arbitrary $m_1,m_2\in(K_9)_1$, which implies that $\varphi$ is trivial \textup{(}since
$(K_9)_1^2=(K_9)_0$\textup{)}.

If $\delta=\frac{1}{2}$ and $m\in(K_9)_1$, then
\[
2\varphi(m)=4\varphi(em)=2e\varphi(m)+\alpha m,
\]
which gives $\varphi(m)=\varphi(m)_0+\alpha m$. Further, we may assume that $\alpha=0$ and
$\varphi(m)=\varphi(m)_0$. Note that $0=\varphi(m^2)=\varphi(m)m$, and hence
\[
\varphi(w)=\alpha_wuw+\beta_wvw,\quad \varphi(z)=\alpha_zuz+\beta_zvz,\quad
\varphi(u)=\alpha_uuz+\beta_uuw,\quad \varphi(v)=\alpha_vvw+\beta_vvz.
\]
It is easy to see that
\[
\varphi(uz)=-(\varphi(u)z+u\varphi(z))=\beta_zz-\beta_uu,
\]
whence
\[
\alpha_uuz+\beta_uuw=\varphi(u)=-\varphi((uz)w)=(\beta_zz-\beta_uu)w+uz(\alpha_wuw+\beta_wvw)
=\beta_ze-\beta_uuw+2\beta_we,
\]
i.e., $\alpha_u=\beta_u=0$. We similarly obtain $\alpha_\gamma=\beta_\gamma=0$ for
$\gamma\in\{z,w,u,v\}$. Thus, $\varphi$ is trivial.

If $\varphi$ is an odd $\delta$-superderivation, then it can readily be proved that
\[
\varphi(e)=\delta\varphi(e)e+\delta e\varphi(e)=\delta\varphi(e)\quad\text{and}\quad\varphi(e)=0.
\]
We also see that, if $\delta\neq\frac{1}{2}$, then
\[
\varphi(a)=\varphi(ea)=\frac\delta2\varphi(a)\quad\text{and}\quad\varphi(m)=2\varphi(em)=2\delta\varphi(m),
\]
which shows that $\varphi$ is trivial.

Let $\delta=\frac{1}{2}$, and let
\[
\varphi(q)=\chi_qe+\alpha_quz+\beta_quw+\gamma_qvz+\mu_qvw
\]
for $q\in\{u,z,v,w\}$. Now it can readily be proved that
\[
\varphi(u)=-\varphi([u,z]w)=\varphi(uz)w+(uz)\varphi(w)=-(\varphi(u)z-u\varphi(z))w+(uz)\varphi(w),
\]
\[
\varphi(u)=\varphi([u,w]z)=-(\varphi(uw)z+(uw)\varphi(z))=(\varphi(u)w-u\varphi(w))z-(uw)\varphi(z).
\]
Obviously,
\[
(\varphi(u)w-u\varphi(w))z-(uw)\varphi(z)\in\zeta\mu_uvw+\langle e,uw,uz,vw\rangle,\quad\zeta=\pm1,
\]
and
\[
-(\varphi(u)z-u\varphi(z))w+(uz)\varphi(w)\in\langle e,uw,uz,vw\rangle,
\]
i.e., $\mu_u=0$. We similarly see that $\gamma_u=0$. Carrying out similar considerations for $v,z,w$, we
obtain $\varphi(q)=\chi_qe+y_qq$, $y_q\in(K_9)_1$.

It is easy to see that
\[
0=2\varphi(e)=2\varphi([z,w])=\varphi(z)w-z\varphi(w)=-\chi_zw+(y_zz)w+\chi_wz-z(y_ww),
\]
and also $(y_zz)w-z(y_ww)\in\langle u,v\rangle$, i.e., $\chi_z=\chi_w=0$. Carrying out similar
manipulations, we obtain $\chi_u=\chi_v=0$.

Further, note that
\begin{multline*}
0=\varphi(e)=-\varphi([u,z][v,w])=\varphi(uz)vw+uz\varphi(vw)\\
=\tfrac12\bigl((\varphi(u)z-u\varphi(z))vw+uz(\varphi(v)w-v\varphi(w))\bigr)
=\zeta_1\beta_uw+\zeta_2\gamma_zv+\zeta_3\gamma_vz+\zeta_4\beta_wu,
\end{multline*}
where $\zeta_i$ takes the values $+1$ and $-1$. Therefore, $\beta_u=\beta_w=\gamma_z=\gamma_v=0$. It
remains to note that
\[
0=-\varphi(e)=-\varphi(uv)=\varphi(u)v-u\varphi(v)=(\alpha_uuz)v-u(\mu_vvw),
\]
which gives $\mu_v=\alpha_u=0$. We can similarly show that $\varphi(z)=\varphi(w)=0$. Thus, $\varphi$ is
trivial.
\end{proof}

\begin{lemma}\label{l2}
Let $V_{\frac{1}{2}}(Z,D)$ be a superalgebra over a field of characteristic~$3$, let $Z$ be the even part of
$V_{\frac{1}{2}}(Z,D)$, and let $\psi$ be a derivation of $Z$ connected with $D$ by the relation
\[
\psi(a)D(b)=D(a)\psi(b)\quad\text{for any } a,b\in Z.
\]
In this case, there is a $z\in Z$ such that $D(z)$ is invertible and $\psi=cD$ for some $c\in Z$.
\end{lemma}

\begin{proof}
It is easy to see that $D(a)\psi=\psi(a)D$. We claim that the element $D(z)$ is invertible for some
$z\in Z$, which will imply that $\psi=(D(z)^{-1}\psi(z))D$, i.e., $\psi$ has the desired property.

Note that $D(a^3)=3D(a)a^2=0$. Hence, it can readily be seen that $a^3Z$ is a $D$-invariant ideal in
$Z$. Thus, either $a^3Z=Z$ or $a^3Z=0$. If $a^3Z=Z$ for every nonzero $a\in Z$,
then $Z$ is a field and every $z$ with $D(z)\neq0$ has the required property. Otherwise, there is
an element $a\in Z$, $a\neq0$, such that $a^3=0$. Consider
\[
R=\{a\in Z\mid a^3=0\}\neq0;
\]
this is an ideal of $Z$ \textup{(}$(a+b)^3=a^3+b^3$ in characteristic~$3$\textup{)} which differs
from $Z$. Let $I\lhd Z$ be a proper ideal. Then every element $a\in I$ is not invertible, and
hence $a^3$ is not invertible; moreover, $a^3Z\neq Z$. Since $a^3Z$ is a $D$-invariant ideal, we obtain
$a^3Z=0$, which implies that $a^3=0$ and $I\subseteq R$. Therefore, $R$ is the largest proper
ideal of $Z$. It is easy to see that $D(R)\not\subseteq R$. Otherwise, $R$ would be a $D$-invariant
ideal, and it would coincide with $Z$. Thus, $D(z)\notin R$ for some $z\in R$. Hence $D(z)$ is
invertible, and this automatically proves the lemma.
\end{proof}

\begin{lemma}\label{l3}
Let $\varphi$ be an odd $\frac{1}{2}$-derivation or an odd $\frac{1}{2}$-superderivation of the superalgebra
$V_{\frac{1}{2}}(Z,D)$ over a field of characteristic~$3$, where $\varphi$ is defined by the conditions
\[
\varphi(a)=\psi(a)x\quad\text{and}\quad\varphi(ax)=\mu(a).
\]
Then $\mu(a)=D(\psi(a))+a\mu(e)$ and $\psi=cD$ for some $c\in A$.
\end{lemma}

\begin{proof}
Note that \textup{(}recall that $\frac{1}{2}=2$ in characteristic~$3$\textup{)}
\[
\psi(ab)x=\varphi(ab)=2(\varphi(a)\cdot b+a\cdot\varphi(b))=(\psi(a)b+a\psi(b))x
\]
for arbitrary $a,b\in A$, i.e., $\psi$ is a derivation of $A$. Moreover,
\begin{multline*}
2\mu(ab)=\varphi(a\cdot bx)=2(\varphi(a)\cdot bx+a\cdot\varphi(bx))=2(\psi(a)x\cdot bx+a\mu(b))\\
=2D(\psi(a))b-2\psi(a)D(b)+2a\mu(b),
\end{multline*}
i.e.,
\begin{equation}\label{e1}
\mu(ab)=D(\psi(a))b-\psi(a)D(b)+a\mu(b),
\end{equation}
which implies for $b=e$ that
\begin{equation}\label{e2}
\mu(a)=D(\psi(a))+a\mu(e).
\end{equation}
Thus, when substituting~\eqref{e2} into~\eqref{e1}, we obtain
\[
D(\psi(ab))+ab\mu(e)=D(\psi(a))b-\psi(a)D(b)+aD(\psi(b))+ab\mu(e),
\]
which gives
\[
D(a)\psi(b)=\psi(a)D(b)
\]
by immediate manipulations. Hence, by Lemma~\ref{l2}, we see that $\psi=cD$ for some $c\in A$.
\end{proof}

\begin{lemma}\label{l4}
Let $\varphi$ be a nontrivial $\delta$-derivation of the superalgebra $V_{\frac{1}{2}}(Z,D)$ over a field of
characteristic $p\neq2$. Then $\delta=\frac{1}{2}$, and
\[
\varphi(y)=(1+p(y))z\cdot y\quad\text{for a chosen } z\in A\setminus Fe.
\]
\end{lemma}

\begin{proof}
It is easy to see that
\[
\varphi(e)=2\delta\varphi(e)e=2\delta\varphi(e)_0+\delta\varphi(e)_1,
\]
i.e., either $\delta\neq\frac{1}{2}$ and $\varphi(e)=0$ or $\delta=\frac{1}{2}$ and $\varphi(e)\in A$.

If $mx\in M$, then
\begin{multline*}
\varphi(mx)=2\varphi(e\cdot mx)=2\delta\varphi(e)\cdot mx+2\delta e\cdot\varphi(mx)\\
=\delta(\varphi(e)m)x+2\delta\varphi(mx)_0+\delta\varphi(mx)_1,
\end{multline*}
which implies that $\varphi(M)=0$ for $\delta\neq\frac{1}{2}$ and $\varphi(mx)_1=(\varphi(e)m)x$ for
$\delta=\frac{1}{2}$.

For an arbitrary $a\in A$,
\[
\varphi(a)=\varphi(e\cdot a)=\delta\varphi(e)\cdot a+\delta e\cdot\varphi(a)
=\delta\varphi(e)\cdot a+\delta\varphi(a)_0+\frac\delta2\varphi(a)_1.
\]
This proves that one of the following conditions holds:
\begin{itemize}
\item[1)] $\delta\neq\frac{1}{2},2$ and $\varphi=0$;
\item[2)] $\delta=2$, $p\neq2,3$, and $\varphi(a)\in M$;
\item[3)] $\delta=\frac{1}{2}$, $p\neq2,3$, and $\varphi(a)=\varphi(e)a$;
\item[4)] $\delta=\frac{1}{2}$ and $p=3$.
\end{itemize}

Let us prove that case~2) gives no nontrivial $2$-derivations. Let $\varphi$ be a $2$-derivation and
$\varphi(a)=\varphi_ax$. Then
\[
0=\varphi(a\cdot x)=2\varphi(a)\cdot x=2\varphi_ax\cdot x=2D(\varphi_a),
\]
and hence, using the fact that $D$ annihilates the elements of the form $\alpha e$ only, we obtain
$\varphi_a=\alpha_ae$. Note that
\[
\varphi_{a^2}x=\varphi(a^2)=4\varphi(a)\cdot a=2(\varphi_aa)x,
\]
i.e., $\varphi_a=0$ for $a\in A\setminus Fe$ \textup{(}and then also for $a=e$, by
linearity\textup{)}. This shows that $\varphi$ is trivial.

Consider case~3), i.e., $\delta=\frac{1}{2}$ and $p\neq2,3$. It is clear that the mapping $\psi$ given by the rule
\begin{equation}\label{e3}
\psi(a)=za,\qquad \psi(mx)=(zm)x\qquad\text{for } z,a,m\in A,
\end{equation}
is a $\frac{1}{2}$-derivation. Therefore, the mapping $\chi$ defined by the rule $\chi=\varphi-\psi$
\textup{(}with $z=\varphi(e)$\textup{)} is also a $\frac{1}{2}$-derivation. We claim that $\chi=0$.

It is easy to see that $\chi(A)=0$ and $\chi(mx)=\varphi(mx)_0$. Moreover,
\[
0=2\chi(x\cdot x)=2\chi(x)\cdot x=\chi(x)x,
\]
i.e., $\chi(x)=0$. Hence we obtain $\chi(ax)=2\chi(a\cdot x)=\chi(a)\cdot x+a\cdot\chi(x)=0$.
Therefore, $\chi{}=0$, and $\varphi$ has the desired form.

Note that
\begin{multline*}
\varphi(ax\cdot bx)=\varphi(D(a)b-aD(b))=z(D(a)b-aD(b))=(ax)\cdot((bz)x)+abD(z)\\
=(ax)\cdot\varphi(bx)+abD(z),
\end{multline*}
for $\varphi$ defined by the rule~\eqref{e3}. Therefore, $\varphi$ is a nontrivial $\frac{1}{2}$-derivation
only if $D(z)\neq0$, i.e., if $z\in A\setminus Fe$. It is clear that the mapping defined in the way
thus described is also a $\frac{1}{2}$-derivation over a field of characteristic $p=3$.

Consider the fourth case, i.e., $\delta=\frac{1}{2}$ and $p=3$. Since the mapping given by formulas~\eqref{e3}
is a $\frac{1}{2}$-derivation of the given superalgebra over a field of characteristic~$3$, we may assume that
$\varphi(A)\subseteq M$ and $\varphi(M)\subseteq A$. Write $\varphi(a)=\psi(a)x$ and
$\varphi(ax)=\mu(a)$. By Lemma~\ref{l3}, we see that $\psi=cD$ for some $c\in A$.

For a mapping $\varphi$ given by the rule
\[
\varphi(a)=(cD(a))x,\qquad \varphi(ax)=D(cD(a))+az\qquad\text{for some } c,z\in A,
\]
to be a $\frac{1}{2}$-derivation, it is necessary that the relation
\[
\varphi(ax\cdot bx)=\tfrac12\bigl(\varphi(ax)\cdot bx+ax\cdot\varphi(bx)\bigr)
\]
be valid. Immediate manipulations show that
\begin{multline*}
(cD(D(a)b-aD(b)))x=\varphi(ax\cdot bx)=\tfrac12\bigl(\varphi(ax)\cdot bx+ax\cdot\varphi(bx)\bigr)\\
=((D(cD(a))+az)b)x+(a(D(cD(b))+bz))x,
\end{multline*}
and hence
\[
acD^2(b)=D(c)D(a)b+aD(c)D(b)+2azb,
\]
which gives $z=0$ for $b=a=e$ and $D(c)D(a)=0$ for $b=e$. Since, by Lemma~\ref{l2}, there is a $w\in A$
such that $D(w)$ is invertible, it follows that $D(c)=0$ and $c=\alpha e$ for some $\alpha\in F$.
Correspondingly, we obtain the relation $\alpha aD^2(b)=0$. Suppose that $\alpha\neq0$. Then
$D^2(b)=0$. This implies that $D(b)=\beta_be$, where $\beta_b\in F$. Note that
\[
\beta_{b^2}e=D(b^2)=2D(b)b=2\beta_bb,
\]
i.e., $\beta_b=0$ for $b\notin Fe$, which gives $D=0$ and contradicts the fact that $D$ is a nonzero
mapping. Thus, $\alpha=0$ and $\psi=\mu=0$.
\end{proof}

\begin{lemma}\label{l5}
Let $\varphi$ be a nontrivial odd $\delta$-superderivation of the superalgebra $V_{\frac{1}{2}}(Z,D)$ over an
algebraically closed field of characteristic $p\neq2$. In this case, $\delta=\frac{1}{2}$ and
\begin{itemize}
\item[1)] if $p\neq3$, then $\varphi(A)=0$ and $\varphi(ax)=az$ for some $z\in A\setminus\{0\}$;
\item[2)] if $p=3$, then $\varphi(a)=(\alpha D(a))x$ and $\varphi(ax)=D(\alpha D(a))+az$ for some
$z\in A$ and $\alpha\in F$ such that the elements $z$ and $\alpha$ do not vanish simultaneously.
\end{itemize}
\end{lemma}

\begin{proof}
Clearly, $\varphi(A)\subseteq M$ and $\varphi(M)\subseteq A$. It is easy to see that
\[
\varphi(e)=2\delta\varphi(e)\cdot e=\delta\varphi(e),
\]
i.e., $\varphi(e)=0$. Let $a\in A$. Then

\begin{gather*}
\varphi(a)=\varphi(e\cdot a)=\delta e\varphi(a)=\frac\delta2\varphi(a),\\
\varphi(ax)=2\varphi(e\cdot ax)=2\delta e\cdot\varphi(ax)=2\delta\varphi(ax).
\end{gather*}

Therefore, one of the following conditions holds:
\begin{center}
$(1)$\ $\delta=\frac{1}{2}$, $p=3$;
$(2)$\ $\delta=\frac{1}{2}$, $p\neq3$;
$(3)$\ $\delta=2$, $p\neq3$;
$(4)$\ $\varphi=0$.

\end{center}
The third case gives $\varphi(M)=0$ and $0=\varphi(a\cdot x)=2\varphi(a)\cdot x$; writing
$\varphi(a)=\varphi_ax$ we get $D(\varphi_a)=0$, so $\varphi_a=\alpha_ae$, and
$\alpha_{a^2}e=2\alpha_aa$ \textup{(}from $\varphi(a^2)=4a\cdot\varphi(a)$\textup{)}, whence
$\varphi_a=0$ for $a\notin Fe$; together with $\varphi(e)=0$ this gives
$\varphi(A)=0$.

The second case gives $\varphi(A)=0$ and $\varphi(ax)=\varphi_a$. Then
\[
0=2\varphi(ax\cdot bx)=\varphi_a\cdot bx-ax\cdot\varphi_b,
\]
which yields $\varphi_ab=a\varphi_b$, i.e., $\varphi_a=a\varphi_e$. For a mapping $\varphi$ defined by
the rule $\varphi(A)=0$ and $\varphi(ax)=za$ to be a $\frac{1}{2}$-superderivation, it is necessary that the
condition
\[
\varphi(a\cdot bx)=\tfrac12\varphi((ab)x)=\tfrac12abz=\tfrac12(0\cdot bx+a\cdot(bz))
=\tfrac12(\varphi(a)\cdot bx+a\cdot\varphi(bx))
\]
be valid. This implies that $\varphi$ is a $\frac{1}{2}$-superderivation. Note that
\[
\varphi(x\cdot x)=0\neq-\tfrac12zx=-x\cdot\varphi(x),
\]
i.e., $\varphi$ is a nontrivial $\frac{1}{2}$-superderivation.

Let us consider the first case, i.e., the case of $p=3$ and $\delta=\frac{1}{2}$, in more detail. Write
$\varphi(a)=\psi(a)x$ and $\varphi(ax)=\mu(a)$. By Lemma~\ref{l3}, $\psi=cD$ for some $c\in A$.

For a mapping $\varphi$ given by the rule
\[
\varphi(a)=(cD(a))x,\qquad \varphi(ax)=D(cD(a))+az
\]
for some $c,z\in A$ to be a $\frac{1}{2}$-superderivation, it is necessary that the relation
\[
\varphi(ax\cdot bx)=\tfrac12\bigl(\varphi(ax)\cdot bx-ax\cdot\varphi(bx)\bigr)
\]
be valid. Immediate manipulations show that
\begin{multline*}
(cD(D(a)b-aD(b)))x=\varphi(ax\cdot bx)=\tfrac12\bigl(\varphi(ax)\cdot bx-ax\cdot\varphi(bx)\bigr)\\
=((D(cD(a))+az)b)x-(a(D(cD(b))+bz))x,
\end{multline*}
which yields $D(c)(D(a)b-aD(b))=0$, i.e., for $b=e$, we readily obtain $D(c)D(a)=0$. By
Lemma~\ref{l2}, there is an element $w\in A$ such that $D(w)$ is invertible, which yields $D(c)=0$, and,
by the conditions in the definition of the derivation $D$ for the algebra $A$, we conclude that
$c=\alpha e$. It is easy to prove that the mapping $\varphi$ given by the rules
\[
\varphi(a)=(\alpha D(a))x,\qquad \varphi(ax)=D(\alpha D(a))+az
\]
for $z\in A$ and $\alpha\in F$ is an odd $\frac{1}{2}$-superderivation. It is clear that, if $\alpha\neq0$ and
$a\neq\beta e$, $\beta\in F$, then
\[
\varphi(a\cdot e)-a\cdot\varphi(e)=\varphi(a)\neq0,
\]
i.e., $\varphi$ is a nontrivial $\frac{1}{2}$-superderivation.
\end{proof}

Note that Lemmas~\ref{l4} and~\ref{l5} give examples of new nontrivial $\frac{1}{2}$-derivations and
$\frac{1}{2}$-superderivations that are not operators of right multiplication.

Recall that new examples of nontrivial $\frac{1}{2}$-derivations for simple unital superalgebras of vector type
$J(\Gamma,D)$ constructed from superalgebras with trivial odd part were presented in~\cite{10}. It turned
out that every operator of right multiplication $R_z$, where $D(z)\neq0$, $z\in\Gamma$, is a nontrivial
$\frac{1}{2}$-derivation and, in turn, all nontrivial $\delta$-derivations of these superalgebras are exhausted
by the operators of right multiplication of this form.

\begin{theorem}\label{t6}
Let $J$ be a simple finite-dimensional Jordan superalgebra over an algebraically closed field of
characteristic $p\neq2$, and let $J$ have a nontrivial $\delta$-derivation or $\delta$-superderivation.
Then $p>2$, $\delta=\frac{1}{2}$, and one of the following holds:

\begin{itemize}
\item[(1)] $J=V_{\frac{1}{2}}(Z,D)$;
\item[(2)] $J=J(B(m),D)$;
\item[(3)] $p=3$ and $J\cong J(V,f)$ is the superalgebra of a superform with $V=V_{\bar1}$, $\dim V=2$.
\end{itemize}
\end{theorem}

\begin{proof}
Note that $\delta$-derivations and $\delta$-superderivations of simple finite-dimensional Jordan
superalgebras over an algebraically closed field of characteristic zero are completely described in the
papers by Kaygorodov~\cite{6,8}. It was proved there that one has no nontrivial $\delta$-derivations and
$\delta$-superderivations for these superalgebras. By the paper of Kaygorodov and Zhelyabin~\cite{10},
nontrivial $\delta$-(super)derivations on simple unital finite-dimensional Jordan superalgebras with
nonsemisimple even part over an algebraically closed field of characteristic $p>2$ can exist only if
$\delta=\frac{1}{2}$ and the superalgebra is of vector type, $J(B(m),D)$ \textup{(}see \cite[Theorem~15]{10};
conversely, $R_z$ with $D(z)\neq0$, $z\in B(m)$, is a nontrivial $\frac{1}{2}$-derivation of $J(B(m),D)$\textup{)}.

Let now $J$ be a simple unital finite-dimensional Jordan superalgebra with semisimple even part over an
algebraically closed field of characteristic $p>2$. By the classification of Racine and Zelmanov~\cite{15},
$J$ is either a simple Jordan algebra \textup{(}trivial odd part\textup{)}, or one of the superalgebras
$M^{(+)}_{m,n}$, $Q(n)^{(+)}$ \textup{(}$n\geq2$\textup{)}, $\mathrm{osp}_{m,2n}$, $P(n)$ \textup{(}$n\geq2$\textup{)},
$J(V,f)$, $D_t$, $K_{10}$, or one of the superalgebras $H_3(B(1,2))$, $H_3(B(4,2))$. The superalgebras
$H_3(B(1,2))$ and $H_3(B(4,2))$ have no nontrivial $\delta$-derivations and $\delta$-superderivations
by \cite[Theorem~15]{10}. For the remaining superalgebras the nontrivial
$\delta$-derivations were excluded in \cite[Theorem~2.5, Lemmas~3.1--3.6 and Theorem~4.4]{6}, and the
nontrivial odd $\delta$-superderivations in \cite[Lemma~8]{8}. Although these results are stated in
characteristic zero, their proofs use only the multiplication tables of the superalgebras and the
invertibility of~$2$ \textup{(}for arbitrary sizes $m,n$ and arbitrary $t$\textup{)}, with one
exception: the computation for odd $\delta$-superderivations of $J(V,f)$ in \cite[Lemma~8]{8} divides
by~$3$. We redo it without this division.

Let $J=J(V,f)=F\cdot1+V$, $V=V_{\bar0}+V_{\bar1}$, $V_{\bar1}\neq0$, and let $\varphi$ be a nonzero odd
$\delta$-superderivation, $\delta\neq0,1$. Since $J$ is unital, $\delta=\frac{1}{2}$ and $\varphi(x)=ux$,
where $u=\varphi(1)\in V_{\bar1}$, $u\neq0$. Since $f$ is nondegenerate, there is $w\in V_{\bar1}$ with
$f(u,w)=1$; then $\varphi(w)=uw=1$. For $x\in V_{\bar0}$ we have $xw=f(x,w)\cdot1=0$ and $\varphi(x)=ux=0$,
whence
\[
0=\varphi(xw)=\frac{1}{2}\bigl(\varphi(x)w+x\varphi(w)\bigr)=\frac{1}{2}\,x,
\]
i.e., $V_{\bar0}=0$. For $y\in V_{\bar1}$ we have
\[
f(w,y)\,u=\varphi(wy)=\frac{1}{2}\bigl(\varphi(w)y-w\varphi(y)\bigr)=\frac{1}{2}\bigl(y-f(u,y)\,w\bigr),
\]
i.e., $y=2f(w,y)\,u+f(u,y)\,w$. Hence $V=V_{\bar1}$ is spanned by $u$ and $w$, and $\dim V=2$. Moreover,
since $f$ is skew-symmetric on $V_{\bar1}$, we have $f(u,u)=0$ and $f(w,u)=-1$, so for $y=u$ we get
$u=-2u$, i.e., $3u=0$, whence $p=3$. Thus $J$ is as in~\textup{(3)}. Conversely, let $p=3$ and
$J=F\cdot1+Fv_1+Fv_2$ with $V=V_{\bar1}=Fv_1+Fv_2$, $v_1v_2=-v_2v_1=1$, $v_1^2=v_2^2=0$. For every
$0\neq v\in V$ the map $t\mapsto vt$ is an odd $\frac{1}{2}$-superderivation of $J$, and it does not lie in
the supercentroid: if $v=\alpha v_1+\beta v_2$ with $\alpha\neq0$, then $v\cdot(v_1v_2)=v$, whereas
$(vv_1)v_2=-\beta v_2$; if $\alpha=0$, then $v\cdot(v_2v_1)=-v=-\beta v_2$, whereas $(vv_2)v_1=0$.

By the paper of Zelmanov~\cite{17}, simple nonunital finite-dimensional Jordan superalgebras over an
algebraically closed field of characteristic $p>2$ are exhausted by the superalgebras of the form $K_3$,
$V_{\frac{1}{2}}(Z,D)$, and the superalgebra $K_9$ in the characteristic~$3$. For the Kaplansky superalgebra
$K_3$, the results on the absence of nontrivial $\delta$-derivations and $\delta$-superderivations follow
from the papers of Kaygorodov~\cite{6,8} independently of the characteristic of the field. Lemma~\ref{l1}
shows that the superalgebra $K_9$ has no nontrivial $\delta$-derivations and $\delta$-superderivations.
Lemmas~\ref{l4} and~\ref{l5} give examples of nontrivial $\delta$-derivations and
$\delta$-superderivations of the superalgebra $V_{\frac{1}{2}}(Z,D)$, which can occur in the case of
$\delta=\frac{1}{2}$ only. With regard to the above comments, this completes the proof of the theorem.
\end{proof}

\section{$\delta$-derivations and $\delta$-superderivations of semisimple finite-dimensional Jordan
superalgebras}

In this section, we obtain a complete description of the $\delta$-derivations and
$\delta$-superderivations of semisimple finite-dimensional Jordan superalgebras over an algebraically
closed field of characteristic distinct from~$2$. In particular, we prove the absence of nontrivial $\delta$-derivations and
$\delta$-superderivations on semisimple finite-dimensional Jordan superalgebras over an algebraically
closed field of characteristic zero.

By a superalgebra with a zero nondivisor we mean a superalgebra $A$ for which there is an element $a$
such that the relation $ax=0$ for some element $x\in A$ implies that $x=0$. The unital superalgebras, the
superalgebras $K_9$ and $V_{\frac{1}{2}}(Z,D)$, and others, are examples of superalgebras of this kind.

\begin{lemma}\label{l7}
Let $\varphi$ be a $\delta$-derivation or a $\delta$-superderivation \textup{(}$\delta\neq0$\textup{)}
of a superalgebra $A$, and let
$A=A_1\oplus A_2,$ where the $A_i$ are superalgebras with a homogeneous zero nondivisor. Then
$\varphi(A_i)\subseteq A_i$.
\end{lemma}

\begin{proof}
Let $e_i$ be a homogeneous zero nondivisor of the superalgebra $A_i$, and let $x_j$ be an
arbitrary element of $A_j$, $j\neq i$. Then
\[
0=\varphi(e_ix_j)=\delta(\varphi(e_i)x_j+(-1)^{p(e_i)p(\varphi)}e_i\varphi(x_j))
=\delta\varphi(e_i)_jx_j+(-1)^{p(e_i)p(\varphi)}\delta e_i\varphi(x_j)_i.
\]
This yields $e_i\varphi(x_j)_i=0$ and $\varphi(x_j)_i=0$, which implies that
$\varphi(A_j)\subseteq A_j$.
\end{proof}

Let $J_i$ be a nonunital superalgebra satisfying the following conditions:
\begin{itemize}
\item[i)] the algebra $(J_i)_0$ has the identity element $e_i$;
\item[ii)] $2e_i\cdot z=z$ for any element $z\in(J_i)_1$;
\item[iii)] if $j\in(J_i)_0$ and $j\cdot(J_i)_1=0$, then $j=0$.
\end{itemize}
The superalgebras $K_9$, $K_3$, and $V_{\frac{1}{2}}(Z,D)$ are examples of superalgebras of this kind
\textup{(}for $K_3=Fe+Fz+Fw$ conditions i)--iii) are immediate from $e^2=e$, $e\cdot z=\frac{1}{2}z$,
$e\cdot w=\frac{1}{2}w$; for $V_{\frac{1}{2}}(Z,D)$ condition~iii) follows from $j\cdot x=\frac{1}{2}jx$; for $K_9$
conditions i)--iii) are checked directly using the multiplication table\textup{)}.

\begin{lemma}\label{l8}
The superalgebra
$J=J_1\oplus\dots\oplus J_n+F\cdot1$
has neither nontrivial $\delta$-derivations nor nontrivial $\delta$-superderivations.
\end{lemma}

\begin{proof}
Let $\varphi$ be a nontrivial $\delta$-derivation or $\delta$-superderivation. By~\cite{6}, a unital
superalgebra can have a nontrivial $\delta$-derivation or $\delta$-superderivation only if $\delta=\frac{1}{2}$.

Since $J$ is unital and $\delta=\frac{1}{2}$, we have $\varphi(y)=\varphi(1)\cdot y$ for all $y\in J$.
Write \begin{center}
$\varphi(1)=\alpha\cdot1+\sum_i(j_0^{(i)}+j_1^{(i)})$ with $j_k^{(i)}\in(J_i)_k$.
\end{center} Since
$J_kJ_i=0$ for $k\neq i$, we get $\varphi(e_i)=\varphi(1)\cdot e_i=\alpha e_i+j_0^{(i)}+\frac{1}{2}j_1^{(i)}$;
on the other hand, $\varphi(e_i)=\varphi(e_i\cdot e_i)=\varphi(e_i)\cdot e_i
=\alpha e_i+j_0^{(i)}+\frac14j_1^{(i)}$. Hence $j_1^{(i)}=0$. For $x\in(J_i)_1$ we have
$\varphi(x)=\varphi(1)\cdot x=\alpha x+j_0^{(i)}\cdot x$ and
\begin{center}$\varphi(x)=2\varphi(e_i\cdot x)=\varphi(e_i)\cdot x+e_i\cdot\varphi(x)
=\frac{1}{2}\alpha x+j_0^{(i)}\cdot x+\frac{1}{2}\alpha x+\frac{1}{2}j_0^{(i)}\cdot x$,\end{center} whence
$j_0^{(i)}\cdot x=0$. By condition~iii), $j_0^{(i)}=0$. Thus $\varphi(1)=\alpha\cdot1$ and
$\varphi=\alpha\,\mathrm{id}$ is trivial.
\end{proof}

By the paper of Zelmanov~\cite{17}, if $J$ is a semisimple finite-dimensional Jordan superalgebra over a
field $F$ of characteristic $p\neq2$, then
\[
J=\bigoplus_{i=1}^s\bigl(J_{i1}\oplus\dots\oplus J_{ir_i}+K_i\cdot1\bigr)\oplus J_1\oplus\dots\oplus J_t,
\]
where $J_1,\dots,J_t$ are simple Jordan superalgebras, $K_1,\dots,K_s$ are extensions of the field $F$,
and $J_{i1},\dots,J_{ir_i}$ are simple nonunital Jordan superalgebras over the field $K_i$,
respectively. Using this result, Lemmas~\ref{l7} and~\ref{l8}, and the papers~\cite{6} and~\cite{8},
where the absence of nontrivial $\delta$-derivations and $\delta$-superderivations of simple
finite-dimensional Jordan superalgebras over an algebraically closed field of characteristic zero was
proved, we can obtain the following theorem.

\begin{theorem}\label{t9}
A semisimple finite-dimensional Jordan superalgebra over an algebraically closed field of characteristic
zero has no nontrivial $\delta$-derivations and $\delta$-superderivations.
\end{theorem}

Using the results of the paper of Zelmanov~\cite{17}, Lemmas~\ref{l7} and~\ref{l8}, and
Theorem~\ref{t6}, we can obtain the following theorem.

\begin{theorem}\label{t10}
Let a semisimple finite-dimensional Jordan superalgebra $J$ over an algebraically closed field of
characteristic $p>2$ have a nontrivial $\delta$-derivation or $\delta$-superderivation. Then
$\delta=\frac{1}{2}$ and $J=J^*\oplus J_*$, where $J_*$ is a simple superalgebra of one of the
types~\textup{(1)}--\textup{(3)} of Theorem~\ref{t6}.

\end{theorem}

\section*{Appendix A. Differences from the published version}
\addcontentsline{toc}{section}{Appendix A. Differences from the published version}

\emph{Theorems~\ref{t6} and~\ref{t10}: statements changed, proof of Theorem~\ref{t6} completed.} The
superalgebra $J(V,f)$ with $V=V_{\bar1}$, $\dim V=2$, in characteristic~$3$   has been added to the list.  

\noindent 
\emph{Lemmas~\ref{l7} and~\ref{l8}}: statement changed.  
\end{document}